\documentclass{amsart}
\usepackage{amsmath}
\usepackage{amsfonts}
\usepackage[latin1]{inputenc}

\newtheorem{theorem}{Theorem}
\newtheorem{corollary}[theorem]{Corollary}
\newtheorem{lemma}[theorem]{Lemma}
\newtheorem{proposition}[theorem]{Proposition}
\newtheorem{definition}[theorem]{Definition}
\newtheorem{conjecture}[theorem]{Conjecture}
\newtheorem{example}[theorem]{Example}
\def\KX{F\langle X \rangle}

\newcommand{\F}{F\langle X|G\rangle}

\newcommand{\Z}{\mathbb{Z}}
\newcommand{\N}{\mathbb{N}}
\newcommand{\zn}{\mathbb{Z}_n}

\newcommand{\U}{\mathcal{U}}

\begin{document}

\title[]{On graded identities of block-triangular matrices with the grading of Di Vincenzo-Vasilovsky}
\author{Lucio Centrone}\address{IMECC, Universidade Estadual de
Campinas, Rua S\'ergio Buarque de Holanda, 651, Cidade Universit\'aria ``Zeferino Vaz'', Distr. Bar\~ao Geraldo, Campinas, S\~ao Paulo, Brazil, CEP
13083-859}\email{centrone@ime.unicamp.br}
\author{Thiago Castilho de Mello}\thanks{T. C. de Mello is supported by the FAPESP grant (No. 2012/16838-0)}
\address{Instituto de Matematica e Estat\'{i}stica, Universidade de S\~{a}o Paulo, Rua do Matão 1010, Caixa Postal 66281, S\~{a}o Paulo, SP, 05315-970
Brasil}\email{tcmello@ime.usp.br}\keywords{}\subjclass[2000]{}
\begin{abstract}
The algebra of $n\times n$ matrices over a field $F$ has a natural $\mathbb{Z}_n$-grading. Its graded identities have been described by Vasilovsky who
extended a previous work of Di Vincenzo for the algebra of $2\times 2$ matrices. In this paper we study the graded identities of block-triangular matrices
with the grading inherited by the grading of $M_n(F)$. We show that its graded identities follow from the graded identities of $M_n(F)$ and from its
monomial identities of degree up to $2n-2$. In the case of blocks of sizes $n-1$ and 1, we give a complete description of its monomial identities, and
exhibit a minimal basis for its $T_{\mathbb{Z}_n}$-ideal.
\end{abstract}\maketitle

\section{Introduction}
Let $F$ be a field and $A$ be an associative algebra over $F$. We say that $A$ is a polynomial identity algebra (PI-algebra, for short), if there exists a
non-zero polynomial $f(x_1,\dots,x_k)$ in the non-commutative variables of the set $X=\{x_1,x_2,\dots\}$, such that $f$ vanishes under any substitution by
elements of $A$. In this case, we say that $f=0$ is a polynomial identity for $A$. If $A$ is a finite dimensional algebra of dimension $d<n$, then it
satisfies the standard polynomial identity of degree $n$
\[s_n(x_1,\dots,x_n)=\sum_{\sigma\in S_n}(-1)^{\sigma}x_{\sigma(1)}\cdots x_{\sigma(n)}=0.\]
In particular, the $n\times n$ matrix algebra over the field $F$, $M_n(F)$, is PI, since it satisfies the standard polynomial $s_{n^2+1}$. In fact by the
well known theorem of Amitsur-Levitzky, $M_n(F)$ satisfies the identity $s_{2n}$.

The set of all polynomial identities of a given algebra $A$, denoted by $T(A)$, is an ideal of $\KX$, the free associative algebra freely generated by $X$,
i.e., the algebra of non-commutative polynomials in the variables of $X$. Moreover, $T(A)$ is invariant under the image of any endomorphism of $\KX$. An
ideal with this property is called a verbal ideal or a T-ideal. Given a set of polynomials $\mathcal{F}\subseteq \KX$ we say that $I$ is the T-ideal
generated by $\mathcal{F}$, if $I$ is the smallest T-ideal containing $\mathcal{F}$. In this case we say that $\mathcal{F}$ is a basis for $I$ or that the
elements of $I$ follow from the elements of $\mathcal{F}$.

One of the central problems in the PI-theory, known as the Specht's problem, is to determine if there exists a finite basis for the polynomial identities
of a given algebra. In charateristic zero, this problem has been completely solved by Kemer \cite{kem1}. Although the existence of a finite basis is known,
it is a very difficult problem to find a finite basis for the T-ideal of a given algebra. If $F$ has characteristic zero, the polynomials $[[x,y]^2,z]$ and
$s_{2n}(x_1,x_2,x_3,x_4)$ form a basis for the polynomial identities of $M_2(F)$. Nevertheless the problem of finding a finite basis for $T(M_n(F))$, for
$n>2$, is still far to be solved. In light of this it is more useful to study ``weaker'' polynomial identities such as identities with trace, identities
with involution, identities with group actions, and graded identities. The graded identities, for example, played an important role in the solution of the
Specht's problem in the work of Kemer. We recall that the Specht's problem has been solved positively in the graded case by Aljadeff and Kanel-Belov
\cite{alk1}.

In this paper we study the graded identities of block-triangular matrix algebras. If $d_1,d_2,\dots,d_m$ are positive integers, such that $d_1+\cdots+
d_m=n$, we denote by $BT(d_1,\dots,d_m;F)$ the subalgebra of the matrix algebra $M_{n}(F)$ consisting of matrices of the type
\[\left(
               \begin{array}{cccc}
                 A_{11} & A_{12} & \cdots & A_{1m} \\
                 0 & A_{22} & \cdots & A_{2m} \\
                 \vdots & \vdots & \ddots & \vdots \\
                 0 & 0 & \cdots & A_{mm} \\
               \end{array}
             \right)\]
where $A_{ij}\in M_{d_i\times d_j}(F)$ for each $i,j$.

The algebra $BT(d_1,\dots,d_m;F)$ is called the algebra of block-triangular matrices of size $d_1$, \dots, $d_m$ over $F$. In a similar way, one defines
$BT(d_1,\dots,d_m;R)$, of block-triangular matrices with entries from an arbitrary algebra $R$.

The block-triangular matrices are quite important in the description of the so called \textit{minimal varieties} of algebras with respect to their
codimension growth, see \cite{giz1, giz2, giz3}. Based on a theorem of Lewin \cite{lew1}, Giambruno and Zaicev \cite{Gia1} have shown that such algebras
have the ``factoring property'', i.e., their T-ideals factor as the product of the T-ideals of their diagonal blocks.

In a fortcoming paper \cite{lctm} the authors give sufficient conditions on the graded PI-algebra $R$ such that $BT(d_1,\dots,d_m;R)$ has the factoring
property. In \cite{dil1} Di Vincenzo and La Scala studied the graded identities of the algebra of block-triangular matrices over a field, graded by a
finite abelian group $G$. They showed that if $R$ is block-triangular with diagonal blocks $A$ and $B$, then $R$ has the factoring property if and only if
one between $A$ and $B$ is $G$-regular. Nevertheless, the $G$-regularity condition does not cover the following case. Let us consider $M_n(F)$ graded by
the Di Vincenzo-Vasilovsky grading and consider the subalgebra of block-triangular matrices, $BT(d_1,\dots,d_m;F)$. In this case $BT(d_1,\dots,d_m;F)$
satisfies monomial identities and does not have the factoring property. Nevertheless, we are still interested in finding a description of monomial
identities and a finite basis for the ideal of graded identities of these algebras. Let us consider a block-triangular matrix algebra such that all its
diagonal blocks are isomorphic to the base field $F$, the so called upper-triangular matrix algebra. In \cite{divkoshval, divkoshval2, koshval} we have a
complete description of the $G$-graded polynomial identities for the upper-triangular matrix algebra when the grading group $G$ is finite abelian.

In this paper we show that the graded polynomial identities of a subalgebra of $M_n(F)$ endowed with the grading of Di Vincenzo-Vasilovsky and generated by
elementary matrices, follow from the Vasilovsky identities of $M_n(F)$ and from its graded monomial identities. In particular, we show that the monomial
identities of $BT(d_1,\dots,d_m;F)$ follow from monomial identities of degree up to $2n-2$. We also give a complete description of the monomial identities
of $BT(n-1,1;F)$.

\section{Graded Structures and Graded Identities of $M_n(F)$}

All fields we refer to are assumed to be of characteristic zero and all algebras we consider are associative and unitary.

Let $(G,\cdot)=\{g_1,\ldots,g_r\}$ be any finite group, and let $F$ be a field. If $A$ is an $F$-algebra, we say that $A$ is a $G$\textit{-graded} algebra
(or $G$-algebra) if there are subspaces $A^g$ for each $g\in G$ such that \[A=\bigoplus_{g\in G}A^g \ \textrm{and} \ A^gA^h\subseteq A^{gh}.\] If $0\neq
a\in A^g$ we say that $a$ is \emph{homogeneous of $G$-degree g} or simply that $a$ has degree $g$, and we write $\deg(a)=g$.

We define a free object; let $\{X^{g}\mid g \in G\}$ be a family of disjoint countable sets of indeterminates. We consider $X=\bigcup_{g\in G}X^{g}$ and we
denote by $\F$ the free associative algebra freely generated by $X$. An indeterminate $x\in X$ is said to be of \emph{homogeneous $G$-degree} $g$, written
$\deg(x)=g$, if $x\in X^{g}$. The homogeneous $G$-degree of a monomial $m=x_{i_1}x_{i_2}\cdots x_{i_k}$ is defined as
$\deg(m)=\deg(x_{i_1})\cdot\deg(x_{i_2})\cdots\deg(x_{i_k})$. For every $g \in G$, we denote by $\F^g$ the subspace of $\F$ spanned by all the monomials
having homogeneous $G$-degree $g$. Notice that $\F^g\F^{g'}\subseteq \F^{gg'}$ for all $g,g' \in G$. Thus \[\F=\bigoplus_{g\in G}\F^g\] is a $G$-graded
algebra, called the \emph{the free associative $G$-graded algebra}. We shall call \textit{$G$-graded polynomials} or, simply, \textit{graded polynomials}
the elements of $\F$. An ideal $I$ of $\F$ is said to be a \textit{$T_{G}$-ideal} if it is invariant under all $F$-endomorphisms $\varphi:\F\rightarrow\F$
such that $\varphi\left(\F^g\right)\subseteq\F^g$ for all $g\in G$. If $A$ is a $G$-graded algebra, a $G$-graded polynomial $f(x_1,\dots,x_t)$ is said to
be a \emph{graded polynomial identity} of $A$ if $f(a_1,a_2,\cdots,a_t)=0$ for all $a_1,a_2,\cdots,a_t\in\bigcup_{g\in G}A^g$ such that $a_k\in
A^{\deg(x_k)}$, $k=1,\cdots,t$. If $A$ has a non-zero graded polynomial identity, we say that $A$ is a \emph{G-graded polynomial identity algebra}
(GPI-algebra). We denote by $T_G(A)$ the ideal of all graded polynomial identities of $A$. We recall that $T_G(A)$ is a $T_G$-ideal of $\F$. If $A$ is
ungraded, i.e., graded by the trivial group, we talk about polynomial identities and T-ideal of $A$. We recall that if two $GPI$-algebras $A$ and $B$
satisfy the same graded identities, i.e., $T_G(A)=T_G(B)$, then they satisfy the same identities, i.e., $T(A)=T(B)$.

Given a subset $Y\subseteq X$ one can talk about the least $T_G$-ideal of $\F$ containing the set $Y$. Such
$T_G$-ideal will be denoted by $\langle Y\rangle ^{T_G}$ and will be called the $T_G$-ideal generated by $Y$. We say that an element of $\langle Y\rangle
^{T_G}$ is a \textit{consequence} of $Y$, or simply, it follows from $Y$. As in the ordinary case, given a $G$-algebra $A$ one of the main problems
in PI-theory is to find a finite set $Y$ such that $T_g(A)=\langle Y\rangle ^{T_G}$. Such a $Y$ is called a basis for the $G$-graded polynomial identities
of $A$. We denote by $\U^G(A)$ the factor algebra $\dfrac{\F}{\F\cap T_G(A)}$ and we shall call it the \textit{relatively free $G$-algebra of $A$}.

We consider now a particular $\Z_n$-grading over $M_n(F)$, the algebra of $n\times n$ matrices with entries from $F$. For each $i, j \in \{1,\dots,n\}$,
let $e_{ij}$ denote the $n\times n$ matrix whose its $(i,j)$-entry is 1 and all the other entries are 0. These matrices, called \textit{elementary
matrices}, form a basis of $M_n(F)$ as a vector space. For each $t\in \mathbb{Z}_n$, let $M_n^{(t)}$ be the subspace of $M_n:=M_n(F)$, spanned by the
elementary matrices $e_{ij}$, such that $j-i\equiv t \; (\text{mod } n) $. The component $0$ consists of diagonal matrices and the component
$t\in\{1,\dots, n-1\}$ consists of matrices of the form
\[\left(
    \begin{array}{ccccccc}
      0           & \cdots & 0       & a_{1,t+1} & 0         & \cdots & 0            \\
      0           & \cdots & 0       & 0       & a_{2,t+2} & \cdots & 0            \\
      \vdots      & \ddots & \vdots  & \vdots  & \vdots    & \ddots & \vdots       \\
      0           & \cdots & 0       & 0       & 0         & \cdots & a_{n-t,n}    \\
      a_{n-t+1,1} & \cdots & 0       & 0       & 0         & \cdots & 0            \\
      0           & \ddots & 0       & 0       & 0         & \cdots & 0            \\
      0           & \cdots & a_{n,t} & 0       & 0         & \cdots & 0            \\
    \end{array}
  \right),
\]where the $a_{i,i+t}$ are elements of $F$. Then $M_n(F)$ is the direct sum
\[M_n(F)=\bigoplus_{t\in \mathbb{Z}_n} M_n^{(t)}.\]
Since the elementary matrices satisfy
\[e_{ij}e_{kl}=\left\{
                 \begin{array}{c}
                   e_{il} \; \text{ if } j=k \\
                    0 \;\;\; \text{ if } j\neq k \\
                 \end{array}
               \right.,
\]
we have $M_n^{(t)}M_n^{(s)}\subseteq M_n^{(t+s)}$ for each $t,s\in \mathbb{Z}_n$, which means that the decomposition $M_n(F)=\bigoplus_{t\in \mathbb{Z}_n}
M_n^{(t)}$ is a $\mathbb{Z}_n$-grading for $M_n(F)$. The graded identities for $M_n(F)$ have been described by Di Vincenzo for $n=2$ (see \cite{div1}) and
by Vasilovsky for arbtrary $n$ (see \cite{vas1}). The main result of Vasilovsky in \cite{vas1} is the following.

\begin{theorem}
All graded polynomial identities of the $\mathbb{Z}_n$-graded algebra $M_n(F)$ follow from the graded identities
\begin{equation}\label{1}
[x_1,x_2]=0, \quad \text{if }\deg(x_1)=\deg(x_2)=0
\end{equation}

\begin{equation}\label{2}
x_1x_2x_3-x_3x_2x_1=0, \quad \text{if } \deg(x_1)=-\deg(x_2)=\deg(x_3)
\end{equation}
\end{theorem}

We denote by $I_n$ the T$_{\Z_n}$-ideal generated by the polynomials from (1) and (2). Let $m=x_1\cdots x_r$ be a monomial of $\F$, then by a
\textit{standard substitution} of $m=x_1\dots x_r$ we mean a substitution $S$ of the form \[x_1\mapsto e_{i_1,j_1}, x_2\mapsto e_{i_2,j_2},\dots,
x_r\mapsto e_{i_r,j_r},\] where $\overline{j_s-i_s}=\deg(x_s)$, so that $e_{i_s,j_s}\in M_n^{\deg(x_s)}$, $s\in \{1,\dots, r\}$ and we denote by $m{|_S}$
the evaluation of $m$ by $S$. In what follows we shall use the following notation. We denote by $m_\sigma$ the monomial $\sigma*m$, where $*$ is the
standard left action of the symmetric group $S_r$ on multilinear polynomials of total degree $r$, i.e., $\sigma*(x_1\dots x_r)=x_{\sigma(1)}\dots
x_{\sigma(r)}$.

In the proof of the theorem above, Vasilovsky uses the next lemma.

\begin{lemma}\label{lemma2}
Let $m=x_1\cdots x_k$ be a graded monomial of total degree $k$. If $\sigma$ and $\tau\in S_k$ and there exists a standard substitution $S$ such that
$m_\sigma|_S=m_\tau|_S\neq 0$, then
\[m_\sigma\equiv m_\tau \;(\text{mod } I_n).\]
\end{lemma}

\section{Graded identities for block-triangular subalgebras of $M_n(F)$}

We consider a graded subalgebra of $M_n(F)$ generated by elementary matrices in order to study its graded identities. In what follows we shall say
\textit{monomial identity} for a monomial that is graded polynomial identity.

\begin{theorem}\label{first}
Let $A$ be a graded subalgebra of $M_n(F)$, generated by elementary matrices. Then the graded polynomial identities of $A$ follow from $I_n$ and from the monomial identities of $A$.
\end{theorem}
\proof Let $J_n$ be the $T_{\Z_n}$-ideal generated by the graded monomial identities of $A$ and by $I_n$.
Since the characteristic of the field $F$ is zero, it is enough to prove that the multilinear graded identities of $A$ are in $J_n$. Let $f(x_1,\dots,x_k)$
be a graded multilinear identity of $A$ and let $r$ be the least non-negative integer such that $f$ can be expressed modulo $J_n$ as a linear combination
of $r$ multilinear monomials:
\[f\equiv \sum_{q=1}^r a_{\sigma_q}m_{\sigma_q} \;(\text{mod }J_n),\]
where $0\neq a_{\sigma_q} \in K$, $\sigma_1$,\dots,$\sigma_r\in S_k$. We show that $r=0$.

Suppose $r>0$, then since $m_{\sigma_1}$ is not a graded identity, it follows that there exists a standard substitution $S$ such that $m_{\sigma_1}|_S\neq
0$. Since for each $q\in\{1,\dots,k\}$, $m_{\sigma_q}|_S\in \{e_{i,j}\,|\,i,j\in \{1,\dots,n\}\}\cup\{0\}$ and
\[a_{\sigma_1}m_{\sigma_1}|_S=\sum_{q=2}^r(-a_{\sigma_q})m_{\sigma_q}|_S,\]
there exists $p\in \{2,\dots,k\}$ such that $m_{\sigma_1}|_S=m_{\sigma_p}|_S\neq 0$. By Lemma \ref{lemma2}, $m_{\sigma_1}\equiv m_{\sigma_p}\; (\text{mod }
I_n)$. Since $I_n\subseteq J_n$, we also have $m_{\sigma_1}\equiv m_{\sigma_p}\; (\text{mod } J_n)$, which implies that
\[f\equiv (a_{\sigma_1}+a_{\sigma_p})m_{\sigma_1}+\sum_{q\neq p, q=2}^r a_{\sigma_q}m_{\sigma_q}\; (\text{mod }J_n) ,\]
which contradicts the minimality of $r$ and we are done. \endproof

Let us consider the algebra of block-triangular matrices, $BT(d_1,\ldots,d_r;F)$, with $d_1+\cdots d_m=n$. If $M_n(F)$ is $G$-graded, then
$BT(d_1,\ldots,d_r;F)$ inherits the structure of a $G$-graded subalgebra of $M_n(F)$. In what follows we shall consider $BT(d_1,\dots,d_m;F)$, with the
grading induced by that of Di Vincenzo-Vasilovsky. It is easy to see that $BT(d_1,\dots, d_m;F)$ satisfies monomial identities. For example, if
$\deg(x_1)=\cdots=\deg(x_n)=\overline{1}$, then $x_1x_2\cdots x_n=0$ is an identity for $BT(d_1,\dots,d_m;F)$. By Theorem \ref{first}, in order to find the
identities of $BT(d_1,\dots,d_m;F)$, it is enough to find its monomial identities. We begin with the following result.

\begin{lemma}\label{zero}
Let $m=x_1\cdots x_t$ be a monomial identity. If there exists $i\in \{1,\dots,t\}$ such that $\deg(x_i)=0$, then $m$ is a consequence of the identity
$m'=x_1\cdots x_{i-1}x_{i+1}\cdots x_t$.
\end{lemma}

We consider the following definition.

\begin{definition}\label{fall}
Let $A,B$ be elements of $M_n(F)$. We say that $A$ has an \emph{empty line} if there exists a line such that its entries are all 0, and we denote by $f_A$
the number of empty lines of $A$. Of course $n-f_A$ is the number of non-empty lines. We say that $AB$ has a \emph{fall between $A$ and $B$}, if
$f_{AB}<f_A$. We shall indicate by $F(A,B)$ the number of falls of $AB$ between $A$ and $B$, i.e., the number $F(A,B)=f_{AB}-f_{A}$.
\end{definition}

We have the following easy lemma.

\begin{lemma}\label{positivefall}
Let $A$, $B$ be homogeneous elements of $M_n(F)$. Then $F(A,B)>0$ if and only if there exists $j$ such that the $j$-th column of $A$ is non-zero and the
$j$-th line of $B$ is zero.
\end{lemma}

\begin{proposition}\label{use}
Let $A,B$ be homogeneous elements of 
$M_n(F)$.
Then $F(A,B)\leq f_B$. Moreover, $F(A,B)=f_B$ 
if and only if for every $i$ such that $A$ has the $i$-th line empty, $B$ has the $(i+\deg(A))$-th line non-empty.

\end{proposition}
\proof Since $A$ and $B$ are homogeneous, we have $A=\sum_{i=1}^n a_{i,j}e_{i,j}$, with $j=i+\deg(A)$ and $B=\sum_{i=1}^n b_{i,j}e_{i,j}$, with
$j=i+\deg(B)$. Then $AB=\sum_{i=1}^n a_{i,k}b_{k,j}e_{i,j}$ with $k=i+\deg(A)$ and $j=i+\deg(A)+\deg(B)$. Hence, $f_{AB}\leq f_A+f_B$, which proves the
inequality. Now we observe that we have equality if and only if $a_{i,k}=0$ implies $b_{k,j}\neq 0$, for every $i$, $j$, and $k$ such that $k=i+\deg(A)$,
and $j=i+\deg(A)+\deg(B)$.
\endproof

We consider the relatively-free $\Z_n$-graded algebra of $BT(d_1,\dots,d_m;F)$. Let $n\in \N$ and set $Y=\{y_{ij}^{(r)}\,|\,i,j=1,\ldots,n;r\geq1\}$. We
shall construct $BT(d_1,\dots,d_m;F)$ as a generic algebra, i.e., as a matrix algebra with entries from the commutative polynomial ring $F[Y]$. For each
$i\in \mathbb{Z}_n$, let $B_i$ be the canonical basis of the component $BT(d_1,\dots,d_m;F)^{(i)}$.

For every $r\geq1$, and $i\in \mathbb{Z}_n$, let us consider \[G_r^{(i)}=\sum_{e_{pq}\in B_i}y_{pq}^{(r)}e_{pq}.\]

Then we have the following classical result.

\begin{theorem}\label{second}
Let $n\geq1$, then $\U^{\zn}(BT(d_1,\dots,d_m;F))\cong F\langle G_r^{(i)}|r\geq1\rangle$.
\end{theorem}

In light of Theorem \ref{second} we may generalize the definition \ref{fall} to the relatively-free graded algebra of $BT(d_1,\dots,d_m;F)$. We note that
Proposition \ref{use} and Lemma \ref{positivefall} also remain valid in this case. From now on let us put $G:=F\langle G_r^{(i)}|r\geq1\rangle$. Moreover,
let us denote by $S$ the set of the $G_r^{(i)}$'s generating $G$. Then we have the following.

\begin{proposition}\label{falls}
Let $m\in \F$ be a monomial of total degree $t$ such that \[m=x_1\cdots x_t,\] where for each $i\in\{1,\ldots,t\}$, $\deg(x_i)=g_i$. Let us consider
\[M=G_1\cdots G_t\] a monomial of $G$ of total degree $t$, where for each $i\in\{1,\ldots,t\}$ $\deg(G_i)=g_i$
. Then $m$ is a monomial identity for $BT(d_1,\dots,d_m;F)$ if and only if there exist $M_1,\ldots,M_{l}$ submonomials of $M$, such that \[M=M_1\cdots
M_{l}\] and for any $i\in\{2,\ldots,l\}$, $F(M_1\cdots M_{i-1},M_i)=\alpha_i$ and $\sum_{i=2}^{l}\alpha_i=n-f_{M_1}$.
\end{proposition}

In light of Proposition \ref{falls} and using the same notations, if $M$ is a monomial of $G$ such that $m$ is a monomial identity, we may decompose $M$ as
\begin{equation} M=G_1G_2C_2G_{i_3}\cdots C_{s-1}G_{i_s}\cdots C_{l-1}G_{i_l}C_{l},\end{equation}\label{3}where $M_1=G_1$, $M_2=G_2$,
$M_s=C_{s-1}G_{i_s}$, if $s>2$ and the $C_s$'s are such that $F(G_1G_2\cdots G_{i_s},C_s)=0$, the $G_{i_s}$'s satisfy $F(G_1G_2\cdots
G_{i_s}C_s,G_{i_{s+1}})>0$ and $l\leq n$. By Lemma \ref{zero}, we suppose $\deg(x_i)\neq 0$, for all $i$. In this case, Lemma \ref{positivefall} implies
$F(G_1,G_2)>0$. If $c_s=\deg(C_s)$ and $z_s=\deg(Z_s)$, then $C_s$ can be written as
\[C_s=G_{s_1}\cdots G_{s_d},\] such that $\deg(G_{s_j})=d_{s_j}$ and $\sum_jd_{s_j}=c_{s}$. We have the following easy lemma.

\begin{lemma}\label{useit3}
Let $m$ and $M$ be as above such that \[M=G_1G_2C_2Z_3\cdots C_{i-1}Z_i\cdots C_{l-1}Z_lC_l.\] Let \[M'=Y'_1Y'_2X'_2Y'_3\cdots Y'_lX'_l\in G,\] where the
$Y'_i$'s are elements of $S$ such that $\deg(Y'_i)=z_i$ and $\deg(X'_i)=c_i$. Then for each $i\in \{2,\dots,l\}$, we have
\[F(G_1G_2C_2\cdots C_{i-1},Z_i)\leq F(Y_1'Y_2'X_2'\cdots X_{i-1}',Y_i').\]
\end{lemma}

\begin{corollary}\label{corollary}
Let $m$, $M$ and $M'$ be as above and let $m'$ be the monomial in $\F$ corresponding to $M'$. Then $m$ is an identity if and only if $m'$ is an identity
for $BT(d_1,\dots,d_m;F)$.
\end{corollary}
\proof Of course $m$ is a consequence of $m'$. So if $m'$ is an identity, also is $m$. Conversely, by Lemma \ref{useit3}, if $m$ is an identity, we have
\[n-f_{G_1}=\sum_{i=2}^l F(G_1C_1Z_2C_2\cdots C_{i-1},Z_i)\leq \sum_{i=2}^l F(Y_1'Y_2'X_2'\cdots X_{i-1}',Y_i')\leq n-f_{Y_1}.\]
Since $f_{G_1}=f_{Y_1}$, we obtain equalities which proves that $m'$ is an identity.
\endproof

As a consequence of Corollary \ref{corollary} we have the following.

\begin{proposition}\label{bound}
Every monomial identity for $BT(d_1,\dots,d_m;F)$ follows from one of degree up to $2n-2$.
\end{proposition}

The above proposition shows that the graded polynomial identities of $BT(d_1,\dots,d_m;F)$, follow from $I_n$ and from monomial identities of degree up to
$2n-2$.

\section{Graded identities of $BT(n-1,1;F)$}

We shall consider the algebra $A$ of block-triangular matrices with two blocks of sizes $n-1$ and $1$, then in light of Theorem \ref{first} it is enough to
find its monomial identities in order to describe its $T_{\Z_n}$-ideal. With abuse of notation, we shall write $j=i$ instead of $j\equiv i\; (\text{mod }
n)$. In what follows we use the following. If $m=m_1m_2$ is a monomial, and $S_1$ is a substitution of $m_1$ and $S_2$ a substitution of $m_2$, we denote
by $S_1S_2$ the substitution of $m$ acting as $S_1$ on $m_1$ and as $S_2$ on $m_2$.

\begin{lemma}\label{subst}
Let $r<n$ be an integer and let $m=x_1\cdots x_r$ be a graded monomial. Then the substitution $S$ given by $x_t\mapsto \sum_{i=1}^{n-1}e_{ij}$, where
$j-i=\deg(x_t)$ yields $m|_S=\sum_{i\in I}e_{ij}$, where $I\subseteq \{1,\dots, n-1\}$, $j-i=\deg(m)$, for every $j$ and $|I|\geq n-r$. Moreover, if
$x_1\cdots x_{r-1}$ does not contain submonomial of degree 0, then $|I|=n-r$.

\end{lemma}
\proof If $r=1$ the assertion is obvious. Suppose now $r>1$ and consider the monomial $m'=x_1\cdots x_{r-1}$, then $m=m'x_r$. Consider a substitution $S'$
on $m'$ such that $m'|_{S'}=\sum_{i\in I'}e_{ij}$, where $I'\subseteq \{1,\dots,n-1\}$, $j-i=\deg(m')$ and $|I'|\geq n-r+1$. We consider now the
substitution $S$ on $m$. We obtain $m|_S=m'|_{S'}x_r|_{S_r}$, where $x_r|_{S_r}=\sum_{i=1}^{n-1}e_{ij}$, with $j-i=\deg(x_r)$. Then $m|_S=\sum_{i\in
I}e_{ij}$, with $j-i=\deg(m')+\deg(x_r)=\deg(m)$ and, $|I|\geq n-r$, in light of Proposition \ref{use}.

To the second part we observe that the substitution $S$ above yields $m|_S=\sum_{i\in I}e_{ij}$, where \[I=\{1,\dots n-1\} \backslash
\{n-\deg(x_1),n-\deg(x_1)-\deg(x_2),\dots,n-\deg(x_1)-\cdots-\deg(x_{r-1})\}\] and if $m'$ does not contain submonomial of degree 0, it means that the
numbers $n, n-\deg(x_1),n-\deg(x_1)-\deg(x_2),\dots,n-\deg(x_1)-\cdots-\deg(x_{r-1})$ are all pairwise distinct $(\text{mod } n)$, which implies that
$|I|=n-r$.
\endproof

As a consequence, we have the following corollary.

\begin{corollary}\label{useitagain}
Let $m=x_1\dots x_r$ be a graded monomial. If $r<n$, then $m$ is not a graded monomial identity for $A$.
\end{corollary}

\begin{lemma}\label{useit2}
Let $m=z_1\cdots z_k x_1\cdots x_r z_{k+1} \cdots z_{n-r}$ be a graded monomial of degree $n$ such that $k+r<n$ and $\deg(x_1\cdots x_r)=0$. If $r>0$
then $m$ is not a graded identity for $A$.
\end{lemma}
\proof If $r>0$, then by Corollary \ref{useitagain} $m_1=z_1\cdots z_k$ and $m_2=z_{k+1}\cdots z_{n-r}$ are not identities for $A$. Moreover by Lemma
\ref{subst}, there exist substitutions $S_1$ and $S_2$ such that $m_1|_{S_1}=\sum_{j\in J_1}e_{ij}$ where $j-i=\deg(m_1)$, $m_2|_{S_2}=\sum_{i\in
J_2}e_{ij}$ where $j-i=\deg(m_2)$ and $|J_1|\geq n-k$  and $|J_2|\geq k+r$. Since $J_1$, $J_2\subseteq \{1,\dots,n-1\}$ and $|J_1|+|J_2|\geq n+r$, we have
$|J_1\cap J_2|\geq r$. By Lemma \ref{subst}, there exists a substitution $S_3$ for $x_1\cdots x_r$ such that $x_1\cdots x_r|_{S_3}=\sum_{i\in I}e_{ii}$ and
$|I|\geq n-r$. Since $I\subseteq\{1,\dots,n-1\}$, and $|J_1\cap J_2|+|I|\geq n$, $I\cap J_1\cap J_2\neq\emptyset$. Then the substitution $S=S_1S_3S_2$,
does not vanish on $m$ and we are done.\endproof

\begin{lemma}\label{useit4}
Let $m=x_1\cdots x_n$ be a graded monomial. Then $m$ is a graded identity for $A$ if and only if $\deg(x_i)\neq 0$, for each $i$, and $m'=x_1\cdots
x_{n-1}$ does not contain submonomials of degree $0$, i.e., if and only if for any $1\leq r< s\leq n-1$ we have $\sum_{t=r}^s\deg(x_t)\neq 0$.
\end{lemma}
\proof We observe that $m$ is a graded identity if and only if it vanishes under the substitution $S$ given by $x_t\mapsto \sum_{i=1}^{n-1}e_{ij}$, where
$j-i=\deg(x_t)$. If $m'$ does not contain a submonomial of degree 0, then $m'|_{S}=e_{i,n}$ for some $i\in\{1,\ldots,n-1\}$ by Lemma \ref{subst} and
$m|_{S}=0$. The converse is also true in light of Lemma \ref{useit2}. \endproof

We give now an example of such result. In what follows, the expression $x_i^{(j)}$ means that the variable $x_i$ has degree $j$.

\begin{example}
Let us consider $BT(4,1;F)$, then the monomial \[x_1^{(1)}x_2^{(2)}x_3^{(3)}x_4^{(4)}x_5^{(4)}=x_1^{(1)}mx_4^{(4)}x_5^{(4)},\]where
$m_1=x_2^{(2)}x_3^{(3)}$, is not a monomial identity for $BT(4,1;F)$, since $m_1$ is a submonomial of $x_1^{(1)}x_2^{(2)}x_3^{(3)}x_4^{(4)}$ of homogeneous
degree 0. On the other hand, the monomial \[x_1^{(1)}x_2^{(2)}x_3^{(1)}x_4^{(3)}x_5^{(4)}\] is a monomial identity for $BT(4,1;F)$ because none of the
submonomials of $x_1^{(1)}x_2^{(2)}x_3^{(1)}x_4^{(3)}$ has homogeneous degree 0.
\end{example}

\begin{proposition}\label{useit5}
Let $k\geq n$ and let $m=x_1\cdots x_k$ be a graded monomial. Then $m$ is an identity for $A$ if and only if $m$ is a consequence of a monomial
identity of degree $n$.
\end{proposition}
\proof  Let $m$ be a graded monomial identity for $A$, then we consider its decomposition in $\U^{\zn}(M_n^k)$ as well as in (\ref{3}),
\[M=G_1G_2C_2G_{i_3}C_3\cdots G_{s_i}C_{i}\cdots G_{s_l}C_l.\]
Write $M_1=G_1$, $M_s=G_{i_s}C_s$, $s\in \{2,\dots, l\}$. We consider the monomial $m'=z_1\cdots z_l\in \F$, such that $\deg(z_s)=\deg(M_s)$, for $s\in
\{1,\dots,l\}$.

It turns out that $m'$ has degree $l\leq n$. Since $f_{G_1}=1$, by Proposition \ref{falls}, for any $s\in\{1,\ldots,l-1\}$, one has $F(M_1M_2\cdots
M_s,M_{s+1})=1$. Consequently, $l=n$. Let us prove now that $m'$ is an identity for $A$. By Lemma \ref{useit4} it is enough to prove that $z_1\cdots
z_{n-1}$ has no submonomials of degree 0. If not, let $Y$ be such a monomial, where $Y=m_i\cdots m_j$. In light of Proposition \ref{falls} we have
$F(m_1\cdots m_{i-1},m_i)=1$. If $m_1\cdots {m_{i-1}}=\sum_{s=1}^{l_1}p_{i_sj_s}(y)e_{i_sj_s}$, where $p_{i_sj_s}(y)$ is a monomial in the $y$'s, then
\[m_1\cdots m_{i-1}{m_i}=\sum_{
\begin{array}{cc}
    s=1\\
    s\neq t
\end{array}
}^{l_1}p_{i_sj_s}(y)e_{i_sj_s},\] where $j_t=n$. Since $\deg(m_1\cdots {m_{i-1}}Y)=\deg(m_1\cdots {m_{i-1}})$, we have $m_1\cdots {m_{i-1}}Y=\sum_{s\in
T}p_{i_sj_s}(y)e_{i_sj_s}$, where $T$ is a proper subset of $\{j_1,\ldots,j_{l_1}\}$. Now the proof follows once we observe that
$F(G_1G_2C_2G_{i_3}C_3\cdots G_{s_j}C_{j}, G_{j+1})=1$ but $n\notin T$. \endproof

In light of Proposition \ref{useit5}, we have the following.

\begin{theorem}
The graded polynomial identities of $A=BT(n-1,1;F)$ follow from the graded identities (\ref{1}), (\ref{2}) and from the monomial identities of degree $n$,
i.e., the monomials $m=x_1\cdots x_n$ such that for any $1\leq r<s\leq n-1$ we have $\sum_{t=r}^s\deg(x_t)\neq 0$. Moreover, such identities form a minimal
basis for the graded identities of $A$.
\end{theorem}

\proof The fact that such identities form a basis for the identities of $A$ follows from Theorem \ref{first} and Proposition \ref{useit5}. To show that
such base is minimal, it is enough to observe that a monomial identity of degree $n$ cannot be a consequence of another monomial of degree $n$ modulo
$I_n$, or we would have a submonomial of degre 0, which does not occur in monomial identities. \endproof

\section{Conclusions}

We observe that the bound obtained in Proposition \ref{bound} does not depend on $d_1,\dots,d_m$ and Proposition \ref{useit5} shows that for $BT(n-1,1;F)$
such bound is exactly $n$. In light of this we state the following conjecture.

\begin{conjecture}
Every monomial identity of $BT(d_1,\dots,d_m;F)$ follows from monomial identities of degree up to $n$.
\end{conjecture}

We also have computed the monomial identities for $BT(2,2;F)$ of degree up to 6 (=$2n-2$). In this case, all monomial identities follow from identities of
degree up to 3. Hence the $n$ in the above conjecture is not minimal, in general.

We can notice that in the proof of Proposition \ref{useit5} we merely used the decomposition (3) with $M_s=G_{i_s}C_s$. One question arises naturally: why
did we not choose $M_s=C_{s-1}G_{i_s}$? Let us consider, for example, $BT(4,1;F)$. It is easy to check that the monomial
\[m=x_1^{(1)}x_2^{(2)}x_3^{(3)}x_4^{(4)}x_5^{(4)}x_6^{(3)}x_7^{(1)}=0\] is a graded identity for $BT(4,1;F)$. Decomposing $m$ as in (\ref{3}), we obtain
\[G_1=x_1^{(1)}, G_2=x_2^{(2)}, G_3=x_3^{(3)}, C_3=x_4^{(4)}x_5^{(4)}, G_{i_4}=x_6^{(3)}, G_{i_5}=x_7^{(1)}.\]
Using the second choice of monomials $M_i$, we obtain
\[m=m_1m_2m_3m_4m_5,\] where $m_1=x_1^{(1)}$, $m_2=x_2^{(2)}$, $m_3=x_3^{(3)}$, $m_4=x_4^{(4)}x_5^{(4)}x_6^{(3)}$ and $m_5=x_7^{(1)}$, but the monomial
\[m'=x_1^{(1)}x_2^{(2)}x_3^{(3)}x_4^{(1)}x_5^{(1)}\] is not a graded identity for $BT(4,1;F)$.
\\

Conversely let us consider $BT(4,4;F)$, then the monomial
\[m=x_1^{(1)}x_2^{(1)}x_3^{(1)}x_4^{(1)}x_5^{(7)}x_6^{(7)}x_7^{(7)}x_8^{(1)}x_9^{(1)}x_{10}^{(1)}x_{11}^{(7)}x_{12}^{(1)}x_{13}^{(1)}=0\] is a graded
identity for $BT(4,4;F)$. Decomposing $m$ as in (\ref{3}), and using the first choice of monomials, we obtain
\[m=m_1m_2m_3m_4m_5m_6m_7m_8,\]where $m_1=x_1^{(1)}$, $m_2=x_2^{(1)}$, $m_3=x_3^{(1)}$, $m_4=x_4^{(1)}$,
$m_5=x_5^{(7)}$, $m_6=x_6^{(7)}$, $m_7=x_7^{(7)}x_8^{(1)}x_9^{(1)}x_{10}^{(1)}x_{11}^{(7)}x_{12}^{(1)}$, $m_8=x_{13}^{(1)}$ but the monomial
\[m'=x_1^{(1)}x_2^{(1)}x_3^{(1)}x_4^{(1)}x_5^{(7)}x_6^{(7)}x_7^{(2)}x_8^{(1)}\] is not a graded identity for $BT(4,4;F)$. Hence the decomposition used to
prove Proposition \ref{useit5} does not hold in general for $BT(n-k,k;F)$.


\begin{thebibliography}{99}
\bibitem{alk1} E. Aljadeff, A. Kanel-Belov, Representability and Specht's problem for $G$-graded algebras, Adv. Math. 225(5) (2010), 2391-2428.
\bibitem{lctm} L. Centrone and T.C. de Mello, \emph{A model for the relatively-free graded algebra of block-triangular matrices with entries from a graded
PI-algebra}, arXiv:1211.6310 (2013).
\bibitem{div1} O. M. Di Vincenzo, On the graded identities of $M_{1,1}(E)$, Israel J. Math. 80 (1992), 323-335.
\bibitem{divkoshval} O. M. Di Vincenzo, P. Koshlukov, A. Valenti, Gradings on the algebra of upper triangular matrices and their graded identities, J.
Algebra, 275 (2004), 550--566.
\bibitem{divkoshval2} O. M. Di Vincenzo, P. Koshlukov, A. Valenti, Gradings and graded identities for the upper triangular matrices over an infinite field,
Groups, Rings, and Group Rings, Lect. Notes Pure Appl. Math., 248 (2006), 91-106.
\bibitem{dil1} O. M. Di Vincenzo, R. La Scala, Block-triangular matrix algebras and factorable ideals of graded polynomial identities, J. Algebra 279 (2004), 260-279.
\bibitem{Gia1} A. Giambruno, M. V. Zaicev, \emph{Polynomial Identities and Asymptotic Methods}, Math. Surveys and Monog. 122 (2005), AMS.
\bibitem{giz3} A. Giambruno, M. V. Zaicev, Codimension growth and minimal superalgebras, Trans. Amer. Math. Soc. 355 (2003), 5091-5117.
\bibitem{giz1} A. Giambruno, M. V. Zaicev, Exponential codimension growth of P.I. algebras: an exact estimate, Adv. Math. 142 (1999), 221-243.
\bibitem{giz2} A. Giambruno, M. V. Zaicev, On codimension growth of finitely generated associative algbras, Adv. Math. 140 (1998), 145-155.
\bibitem{kem1} A. R. Kemer, Varieties and $\Z_2$-graded algebras (Russian), Izv. Akad. Nauk SSSR, Ser. Mat. 48 (1984), 1042-1059. Translation: Math. USSR, Izv. 25 (1985), 359-374.
\bibitem{koshval} P. Koshlukov, A. Valenti, Graded identities for the algebra of {$n\times n$} upper triangular matrices over an infinite field, Internat. J. Algebra
Comput. 13 (2003) 517--526.
\bibitem{lew1} J. Lewin, A matrix representation for associative algebras I, Trans. Amer. Math. Soc. 188 (1974), 293-308.
\bibitem{vas1} S. Y. Vasilovsky, $\zn$-graded polynomial identities of the full matrix algebra of order $n$, Proc. Amer. Math. Soc. 127(12) (1999), 3517-3524.
\end{thebibliography}
\end{document}